\documentclass[11pt,oneside,english]{amsart}
\usepackage[T1]{fontenc}
\usepackage[latin9]{inputenc}
\usepackage{textcomp}
\usepackage{amsthm}
\usepackage{esint}

\makeatletter
\numberwithin{equation}{section}
\numberwithin{figure}{section}
\theoremstyle{plain}
\newtheorem{thm}{Theorem}
  \theoremstyle{plain}
  \newtheorem*{conjecture*}{Conjecture}

\def\makebbb#1{
    \expandafter\gdef\csname#1\endcsname{
        \ensuremath{\Bbb{#1}}}
}\makebbb{R}\makebbb{N}\makebbb{Z}\makebbb{C}\makebbb{H}\makebbb{E}\makebbb{H}\makebbb{P}\makebbb{Q}\makebbb{K}\makebbb{E}

\usepackage{babel}

\usepackage{babel}

\makeatother

\usepackage{babel}

\makeatother

\usepackage{babel}

\begin{document}

\title{The projective space has maximal volume among all toric Kähler-Einstein
manifolds}

\author{Robert J. Berman, Bo Berndtsson}

\email{robertb@chalmers.se, bob@chalmers.se}

\curraddr{Mathematical Sciences - Chalmers University of Technology and University
of Gothenburg - SE-412 96 Gothenburg, Sweden }
\begin{abstract}
We prove a conjecture saying that complex projective space has  maximal
volume (degree) among all toric Kähler-Einstein manifolds of dimension
$n.$ The proof is inspired by our recent work on sharp Moser-Trudinger
and Brezis-Merle type inequalities for the complex Monge-Ampère operator,
but is essentially self-contained. 
\end{abstract}
\maketitle
Let $X$ be an $n-$dimensional complex manifold $X$ which is \emph{Fano}
(i.e. its first Chern class $c_{1}(X)$ is ample/positive). For some
time it was expected that top-intersection number $c_{1}(X)^{n},$
also called the \emph{degree} of $X,$ is maximal for the $n-$dimensional
complex projective space, i.e. \begin{equation}
c_{1}(X)^{n}\leq(n+1)^{n},\label{eq:conj ineq}\end{equation}
There are now counterexamples to this bound. For example, as shown
by Debarre (see page 139 in \cite{de}), even in the case when $X$
is \emph{toric} (i.e. it admits an effective holomorphic action of
the complex torus $(\C^{*})^{n}$ with an open dense orbit) there
is no universal polynomial bound in $n$ on the $n-$th root of the
degree $c_{1}(X)^{n}.$ However a specific conjecture concerning the
toric case (see \cite{n-p} and references therein), says that the
bound above holds when $X$ is \emph{Kähler-Einstein }(i.e. it admits
a Kähler metric $\omega$ with constant Ricci curvature). In this
note we will confirm this conjecture:
\begin{thm}
Let $X$ be an $n-$dimensional smooth toric F variety which admits
a Kähler-Einstein metric. Then its first Chern class $c_{1}(X)$ satisfies
the following upper bound \[
c_{1}(X)^{n}\leq(n+1)^{n}\]

\end{thm}
As pointed out in \cite{n-p} one of the motivations for the bound
above in the toric setting is another more general conjecture of Ehrhart
in the realm of convex geometry, which can be seen as a variant of
Minkowski\textquoteright{}s first theorem for non-symmetric convex
bodies:
\begin{conjecture*}
(Ehrhart). Let $P$ be an n-dimensional convex body which contains
precisely one interior lattice point. If the point coincides with
the barycenter of $P$ then \[
\mbox{Vol}(P)\leq(n+1)^{n}/n!\]

\end{conjecture*}
The case when $n=2$ was settled by Ehrhart, as well as the special
case of simplices in arbitrary dimensions \cite{e}. As explained
in the survey \cite{g-w} the best general upper bound in Ehrhart's
conjecture, to this date, is $\mbox{Vol}(P)\leq(n+1)^{n}(1-(n-1)/n)^{n}(\leq(n+1)^{n}).$
According to the well-known dictionary between polarized toric varieties
$(X,L)$ and convex lattice polytopes $P$ \cite{do,de2:} the previous
theorem confirms Ehrhart's conjecture for \emph{lattice polytopes}
$P$ of the form \begin{equation}
P=\{x\in\R^{n}:\,\,\left\langle l_{i},x\right\rangle \leq1\}\label{eq:poly}\end{equation}
where $l_{i}$ are primitive lattice vectors and such that $P$ is
\emph{Delzant,} i.e. any vertex of $P$ meets precisely $n$ facets
and the corresponding $n$ vectors $l_{i}$ generate the lattice $\Z^{n}.$
In any such polytope the origin $0$ is indeed the unique lattice
point and, as shown by Wang-Zhou \cite{w-z}, $0$ is the barycenter
of $P$ precisely when the corresponding toric Fano manifold $X$
admits a Kähler-Einstein metric. The cases up to $n\leq8$ have previously
been confirmed by computer assistance (as announced in \cite{n-p}),
using the classification of Fano polytopes for $n\leq8$ \cite{oe}.

It should also be pointed out that, as shown in \cite{g-m-s-y}, Bishop's
volume estimate for Einstein metrics, applied to the unit circle-bundle
in the canonical line bundle $K_{X}\rightarrow X$ translates, to
the inequality \begin{equation}
c_{1}(X)^{n}\leq(n+1)^{n}\frac{(n+1)}{I(X)}\label{eq:bishop}\end{equation}
where $I(X)$ is the Fano index of $X,$ i.e. the largest positive
integer $I$ such that $c_{1}(X)/I$ is an integral class in the Picard
group of $X.$ As is well-known $I(X)\leq n+1$ with equality precisely
for $X=\P^{n}$ (see for example page 245 in \cite{k} ) and hence
the previous theorem improves on the inequality \ref{eq:bishop} in
the case when $X$ is toric. 

The idea of the proof which is inspired by our previous work \cite{b-b},
is that if $X$ admits a Kähler-Einstein metric $\omega$ violating
the inequality \ref{eq:conj ineq} (where $c_{1}(X)^{n}/n!$ coincides
with the volume of $\omega)$ then we obtain a remarkably good Moser-Trudinger
inequality for $T^{n}-$invariant plurisubharmonic functions on a
sufficiently large domain $\Omega$ in $\C^{n},$ equivariantly embedded
into $X.$ The contradiction is obtained by showing that the Moser-Trudinger
inequality is simply too good to be true. To this end we show that
the Moser-Trudinger inequality implies a lower bound on the integrability
index of $T^{n+1}-$invariant plurisubharmonic functions $u$ on the
product domain $\Omega':=\Omega\times D$ in $\C^{n+1},$ where $D$
is the unit-disc. But the bound is violated by the pluricomplex Green
function of $\Omega'$ with a pole at the origin, which gives the
desired contradiction.

\subsubsection*{Generalizations}

In fact, our arguments show that the theorem above is valid for any
Kähler-Einstein manifold $X$ which is an $S^{1}-$equivariant compactification
of $\C^{n}.$ However, the symmetry under the full torus simplifies
some of the technical aspects of the proof. 

Moreover, our method of proof can also be modified to handle the case
of possibly \emph{singular} toric Fano varieties, i.e. $-K_{X}$ is
an ample $\Q-$Cartier divisor (equivalently: the corresponding polytope
$P$ in \ref{eq:poly} is merely assumed to be rational and not necessarily
Delzant). To make the connection with the Ehrhart conjecture we are
then led to extend the result of Wang-Zhou \cite{w-z} to the case
of a general Fano varieties using the notation of Kähler-Einstein
metrics on Fano manifolds with log-terminal singularities very recently
introduced in \cite{bbegz}. Even more generally we obtain the existence
of Kähler-Einstein metrics on the complex torus $\C^{*n}$ whose boundary
behavior is determined by a given convex body $P.$ This is a consequence
of the following theorem:
\begin{thm}
Let $P$ be a bounded convex body continaining $0$ in its interior.
Then $0$ is the barycenter of $P$ if and only if there is a convex
function $\phi$ on $\R^{n}$ such that \[
MA_{\R}(\phi)=e^{-\phi}dx\]
 and such that the gradient image of $\phi$ is $P.$ Moreover, the
solution is unique modulo the action of the group $\R^{n}$ by translations.
\end{thm}
Here $MA_{\R}(\phi)$ denotes the real Monge-Ampère measure of the
convex function $\phi$ and $dx$ is the usual Euclidean volueme form.
Details will appear elsewhere \cite{b-b2}.

\subsubsection*{Acknowledgments}

We are grateful to Benjamin Nill for helpful comments.

\subsection{Proof of the theorem}

Assume for a contradiction that $c_{1}(X)^{n}>(n+1)^{n}.$ As is well-known
any smooth and compact toric variety may be realized as an equivariant
compactification of $\C^{n}$ with its standard torus action (see
for example page 10 in \cite{do}). In other words, we may embed $F:\,\C^{n}\rightarrow X$
as an open dense set in $X$ in such a way that the action of the
unit-torus $T^{n}$ and its complexification is preserved. Let now
$\omega$ be an $T^{n}-$invariant Kähler-Einstein metric on $X.$
We can then write \[
F^{*}\omega=dd^{c}u:=\frac{i}{2\pi}\partial\bar{\partial}u\]
 for a smooth $T^{n}-$invariant function $u$ on $\C^{n}$ satisfying
the Kähler-Einstein equation \[
(dd^{c}u)^{n}=Ce^{-u}dV\]
 for a positive constant $C.$ We may rewrite $C=V_{X}/\int_{\C^{n}}e^{-u}dV$,
where \[
V_{X}=\int_{\C^{n}}(dd^{c}u)^{n}\]
which coincides with the top-intersection number $c_{1}(X)^{n}.$
Let now $\Omega_{R}$ be the set where $u<R$ and note that the sets
$\Omega_{R}$ exhaust $\C^{n},$ i.e. $u$ is proper. Indeed, by symmetry
$0$ is a critical point for $u$ and since $dd^{c}u>0$ $u$ is strictly
convex in the logarithmic coordinates $p_{i}=\log|z_{i}|^{2}$ it
follows that $u\rightarrow\infty$ as $|p|\rightarrow\infty.$ We
can hence fix $R$ sufficiently large so that \[
V_{\Omega_{R}}:=\int_{\Omega_{R}}(dd^{c}u)^{n}>(n+1)^{n}.\]
 Writing $\Omega:=\Omega_{R}$ and replacing $u$ by $u-R$ we then
obtain a $T^{n}-$invariant smooth plurisubharmonic (psh for short)
function $u$ (i.e. $dd^{c}u\geq0)$ solving the following equation:
\begin{equation}
(dd^{c}u)^{n}=V_{\Omega}\frac{e^{-u}dV}{\int_{\Omega}e^{-u}dV},\,\,\,\, u=0\,\,\mbox{on\,\ensuremath{\partial\Omega}}\label{eq:m-a eq}\end{equation}
on the $T^{n}-$invariant domain $\Omega$ which is \emph{hyperconvex},
i.e. it admits a negative continuous psh exhaustion function (namely
$u).$ We will denote by $\mathcal{H}_{0}(\Omega)$ the space of all
psh functions on $\Omega$ which are continuous up to the boundary,
where they are assumed to vanish. Its $T^{n}-$invariant subspace
will be denoted by $\mathcal{H}_{0}(\Omega)^{T^{n}}.$ 

Using the previous equation the following Moser-Trudinger type inequality
can now be established on $\Omega$ (which has $u$ above as an extremal):
there is a positive constant $C$ such that \[
\mbox{(M-T)\,\,\,\,}\log\int_{\Omega}e^{-u}dV\leq\frac{1}{V_{\Omega}}\frac{1}{(n+1)}\int_{\Omega}(-u)(dd^{c}u)^{n}+C\]
 for any $u\in\mathcal{H}_{0}(\Omega)^{T^{n}}.$ The proof may be
obtained by repeating the proof of Theorem 1.4 in \cite{b-b}, but
for completeness we have given a slightly simplified proof in section
\ref{sub:Proof-of-the} below, which takes advantage of the full $T^{n}-$symmetry
(while the general argument in \cite{b-b} only requires $S^{1}-$symmetry).

Next, from the previous Moser-Trudinger type inequality we deduce
another inequality of Brezis-Merle type on the hyperconvex domain
$\Omega'=\Omega\times D$ in $\C^{n'}$ where $n'=n+1:$ there is
positive constant $A$ such that \[
\mbox{(B-M);\,\,\,\,}\int_{\Omega'}e^{-(V_{\Omega}n')^{1/n'}u}dV\leq A\left(1-\int_{\Omega'}(dd^{c}u)^{n'}\right)^{-1}\]
for any $u\in\mathcal{H}_{0}(\Omega')^{T^{n'}}$ such that $\int_{\Omega'}(dd^{c}u)^{n'}<1.$
In particular, since by assumption $V_{\Omega}>(n+1)^{n},$ this forces
\[
\int_{\Omega'}e^{-n'u}dV\leq A'<\infty\]
 for any $u\in\mathcal{H}_{0}(\Omega')^{T^{n'}}$ such that $\int_{\Omega'}(dd^{c}u)^{n'}=1.$
More generally, taking limits the previous inequality holds for all
$u\in\mathcal{F}(\Omega')^{T^{n'}},$ where $\mathcal{F}(\Omega')$
is Cegrell's class, which by definition consists of all psh functions
$u$ on $\Omega$ which are decreasing limits of elements $u_{j}$
in $\mathcal{H}_{0}(\Omega')$ with a uniform upper bound on the Monge-Ampère
masses $\int_{\Omega}(dd^{c}u)^{n}$ (see \cite{b-b} and references
therein)..

The desired contradiction will now be obtained by exhibiting a function
violating the previous inequality. To this end we simply let $u:=g$
be the pluricomplex Green function for $\Omega'$ with a pole at $0:$
\begin{equation}
g(z):=\sup\left\{ u(z):\,\,\, u\in(PSH\cap\mathcal{C}^{0})(\Omega'-\{0\}):u\leq0,\,\,\, u\leq\log|z|^{2}+O(1)\right\} \label{eq:g as sup}\end{equation}
 As is well-known \cite{kl} $g$ is continuous up to the boundary
on $\Omega'$ a part from a singularity at $z=0$ and satisfies \begin{equation}
(dd^{c}g)^{n}=\delta_{0}\,\mbox{\,\ on\,}\Omega'-\{0\},\,\,\,\, g=\log|z|^{2}+O(1)\label{eq:eq for g}\end{equation}
In particular $(dd^{c}g)^{n}=1$ and $\int_{\Omega'}e^{-n'g}dV=\infty.$
Finally, the proof is concluded by noting that $g$ is $T^{n}-$invariant.
Indeed, since $0$ is invariant under the action of $T^{n}$ it preserves
the convex class of functions where the sup in \ref{eq:g as sup}
is taken and hence the sup $g$ must be $T^{n}-$invariant. Alternatively
one can also invoke the uniqueness of solutions to \ref{eq:eq for g}
(in the class of functions with the same regularity properties as
$g$).

\subsection{\label{sub:Proof-of-the}Proof of the Moser-Trudinger type inequality
(M-T)}

Let \[
\mathcal{G}(u):=\log\int_{\Omega}e^{-u}dV+\frac{1}{V_{\Omega}}\frac{1}{(n+1)}\int_{\Omega}u(dd^{c}u)^{n}\]
whose Euler-Lagrange equation (i.e. the critical point equation $d\mathcal{G}_{|u}=0)$
is precisely the complex Monge-Ampère equation \ref{eq:m-a eq}. Given
$u_{0}$ and $u_{1}$ in $\mathcal{H}_{0}(\Omega)$ there is a unique
\emph{geodesic} $u_{t}$ connecting them in $\mathcal{H}_{0}(\Omega)$
which may be defined as the unique solution to the following Dirichlet
problem for the Monge-Ampère equation: setting $U(z.t):=u_{t}(z),$
where now $t$ has been extended to a complex strip $\mathcal{T}$
by imposing invariance in the imaginary $t-$direction, we have, for
$M:=\Omega\times\mathcal{T},$ that $U\in\mathcal{C}^{0}(\bar{M})\cap PSH(M)$
and \[
(dd^{c}U)^{n+1}=0,\,\,\,\mbox{in\,\ \ensuremath{M}}\]
and on $\partial M$ the function $U$ coincides with the boundary
data determined by $u_{0}.$ Alternatively, $U$ can be directly defined
as the sup over all $V$ in $\mathcal{C}^{0}(\bar{M})\cap PSH(M)$
restricting to the given boundary data on the boundary. In particular,
if $u_{0}$ and $u_{1}$ are in $\mathcal{H}_{0}(\Omega)^{T^{n}}$
then so is $u_{t.}$ In fact, in the $T^{n}-$invariant case $u_{t}$
may alternatively be obtained by Legendre transform considerations
(compare \cite{s-z}). Indeed, letting $p_{i}:=\log|z_{i}|^{2}$ we
can identify $u(z)$ with a convex funciton on a domain in $\R^{n}$
(with coordinates $p)$ that we, abusing notation slightly, write
as $u(p).$ Moreover, we can extend $u(p)$ to become a smooth convex
funtions on all of $\R^{n}$ such that $u(p)=C\max_{i}\{p_{i}\}+O(1)$
as $p\rightarrow\infty$ and by a simple approximation argument we
may as well assume that $u_{i}(p)$ is smooth and strictly convex.
Then $u_{t}(z)$ may be obtained by connecting the Legendre transforms
$u_{i}^{*}$ of $u_{0}(p)$ and $u_{1}(p)$ by an affine curve and
then taking the Legendre transform again, i.e. $u_{t}(x)=(u_{0}^{*}(1-t)+u_{1}^{*}t)^{*}.$
In particular $u_{t}$ is always smooth in the $T^{n}-$invariant
case if its end points are smooth and strictly convex, which simplifies
some of the technical points of the proof in \cite{b-b}.

As shown in \cite{bern1} the functional \[
t\mapsto\log\int_{\Omega}e^{-u_{t}}dV\]
 is concave along any geodesic as long as the domain $\Omega$ is
$S^{1}-$invariant. In the $T^{n}-$invariant case this fact can also
be deduced from the Prekopa-Leindler inequality for convex functions
on $\R^{n}.$ Indeed, since $u(p)$ is convex in $p$ and vanishes
on $\partial\Omega$ we may extend it to a convex function on all
of $\R^{n}$ by letting it be equal to $\infty$ on the complement
of $\Omega.$ Then the required concavity follows from the Prekopa-Leindler
inequality which says that $-\log\int_{\R^{n}}e^{-v_{t}}dp_{1}\wedge\cdots\wedge dp_{n}$
is convex in $t$ if $v$ is convex in $(p,t)$ (take $v_{t}(p)=u_{t}(p)+\sum_{i}p_{i}).$ 

Next, we note that\[
\mathcal{E}(u):=\int_{\Omega}u(dd^{c}u)^{n}\]
 is affine along a geodesic $u_{t}.$ Indeed, letting $t$ be complex
a direct calculation gives \begin{equation}
dd^{c}\mathcal{E}(u_{t})=\int_{\Omega}(dd^{c}u)^{n+1}\label{eq:psh-forward formula}\end{equation}
which, by definition, vanishing if $u_{t}$ is a geodesic. All in
all this means that $\mathcal{G}(u_{t})$ is concave along a geodesic.
Letting now $u$ be an arbitrary element in $\mathcal{H}_{0}(\Omega)^{T^{n}}$we
take $u_{t}$ to be the geodesic connecting the solution $u_{0}$
of equation \ref{eq:m-a eq} (obtained from the Kähler-Einstein metric
on $X)$ and $u_{1}=u.$ Then $\mathcal{G}(u_{t})$ has a critical
point at $t=0,$ i.e. its right derivative vanishes for $t=0$ and
hence by concavity $\mathcal{G}(u_{1})\leq\mathcal{G}(u_{0})$ which
concludes the proof of the M-T inequality with $C=\mathcal{G}(u_{0}).$

\subsection{Proof of the Brezis-Merle type inequality (B-M)}

Let $u(z,t)$ be an element in $\mathcal{H}_{0}(\Omega\times D)^{T^{n+1}}.$
Applying the M-T inequality established above on $\Omega$ for $t$
fixed gives \[
\int_{\Omega}e^{-u_{t}}dV\leq\exp(-\frac{1}{V_{\Omega}}\frac{1}{(n+1)}\mathcal{E}(u_{t}))\]
By \ref{eq:psh-forward formula} $\mathcal{E}(u_{t})$ is psh for
$t\in D,$ vanishing on the boundary and hence applying the one-dimensional
Brezis-Merle inequality on the unit-disc: \[
\int_{D}e^{-v}dV\leq A(1-\int_{\Omega}dd^{c}v)^{-1},\]
for $v\in\mathcal{H}_{0}(D)$ such that $\int_{\Omega}dd^{c}v<1,$
gives \[
\int_{\Omega\times D}e^{-u}dV\leq A(1-\frac{1}{V_{\Omega}}\int_{\Omega\times D}(dd^{c}u)^{n+1})^{-1}\]
so that rescaling $u$ concludes the proof. Note that the one-dimensional
Brezis-Merle inequality used above is a simple consequence of Green's
formula. In fact, we only need the basic fact that $\int_{D}e^{.-v}<\infty$
if $v$ is subharmonic and bounded on $D$ and $\int_{D}dd^{c}v<1.$

\end{document}